\documentclass[11pt]{article}

\usepackage{amssymb,amsmath,latexsym}
\usepackage[dvips]{graphics}
\usepackage{mathrsfs}

\newtheorem{thm}{Theorem}[section]

\newtheorem{lemma}[thm]{Lemma}
\numberwithin{equation}{section}
\renewcommand{\baselinestretch}{1.4}
\def\ZZ{\mathbb{Z}}

\def\HH{\mathbb{H}}

\def\EE{\mathbb{E}}

\def\VV{\mathscr{V}}

\def\RR{{\mathbb{R}}}

\def\<{\langle}
\def\>{\rangle}
\def\pf{\noindent{\bf Proof.} }

\def\VV{\mathbb{V}}

\def\qed{{\hfill $\Box$\medskip}}

\allowdisplaybreaks[4]

\input{epsf.sty}
\usepackage{epsfig}

\begin{document}
\title{\bf A study  of random walks on  wedges}
\author{Xinxing Chen \footnote{Department of Mathematics, Shanghai Jiaotong University, Shanghai 200240, China;
Research  partially supported  by the NSFC grant No. 11001173.
chenxinx@sjtu.edu.cn}
 } \maketitle

\begin{abstract}
 In this paper we  develop the idea of
Lyons and  gives a simple criterion for the recurrence and the
transience.
  We also show that a wedge has the infinite collision property
if and only if it is a recurrent graph.

 \end{abstract} \noindent{\bf 2000 MR subject
classification:} 60K


\section{Introduction}
Let us recall briefly the definition of a wedge of $\ZZ^{d+1}$. Let
$f_1,\cdots, f_d$ be a collection of $d$ increasing functions from
$\ZZ^+$ to $\RR^+\cup\{+\infty\}$. They induces a wedge,
Wedge$(f_1,\cdots, f_d)=(\VV,\EE)$, which has vertex set
$$
\VV=\{(x_1,\cdots,x_d, n)\in \ZZ^{d+1}: n\ge 0,~ 0\le x_i\le
f_{i}(n) {\rm ~for~ each~}1\le i\le d\}
$$
and edge set
$$
\EE=\{[u,v]: \|u-v\|_1=1, u,v\in \VV\}.
$$
Is a wedge recurrent or transient? (A locally finite connect graph
is called transient or recurrent according to the type of simple
random walk on it.) Lyons\cite{TL} first give the result that
suppose  (A) holds, then Wedge($f_1,\cdots,f_d$) is recurrent if and
only if
\begin{equation}\label{e:1.0}
\sum_{n=0}^\infty\prod_{i=1}^d\frac{1}{f_{i}(n)+1}=\infty.
\end{equation}
Where \\
$~~~~~~~~~~~~$(A):
 $f_i(n+1)-f_i(n)\in \{0,1\}$  for all $1\le i\le d$
and all $n\ge 0.$\\
Readers can refer to \cite{AB}\cite{Peres} for more background about
wedge and the reference therein.

 We develop the idea of Lyons in this paper. However, our result does not rely on the  condition (A).
Define $d$ increasing integer valued functions $h_1,\cdots, h_d$.
Let $ h_i(0)=0$ for each $1\le i\le d$. For each $1\le i\le d$ and
$n\ge 1$, if $h_i(n-1)+1>f_i(n)$ then let
$$
h_i(n)=h_{i}(n-1);
$$
otherwise, if $h_i(n-1)+1\le f_i(n)$ then let
$$
h_i(n)=h_i(n-1)+1.
$$
Then we have our first result.
\begin{thm}\label{t:1.2}
Wedge$(f_1,\cdots,f_d)$ is recurrent if and only if
\begin{equation}\label{e:1.1}
\sum_{n=0}^\infty\prod_{i=1}^d\frac{1}{h_{i}(n)+1}=\infty.
\end{equation}
\end{thm}
{\it Example.} Suppose $d=2$, $f_1(x)=2^x$ and $f_2(x)=\log(x+1)$.
Obviously (\ref{e:1.0}) does not succeed. On the other hand,
$h_1(n)=n$ and $h_2(n)=[\log(n+1)]$. Then (\ref{e:1.1}) holds and
Wedge($f_1,f_2$)
is recurrent.\\

Now we turn to another question. As usual, we say that a graph has
the infinite collision property if
 two independent simple random walks on the graph will  collide infinitely many
 times, almost surely.   Likewise
we say that a graph has the finite collision property if  two
independent simple random walks on the graph collide finitely many
times almost surely. It is interesting to known whether or not  a
graph has the infinite collision property. Refer to
Polya\cite{Polya}, Liggett\cite{LT} and Krishnapur \& Peres\cite{KP}
for details. To my interest is the type of a wedge. Other graphes,
such as wedge combs, trees or random environment, are studied in
\cite{BP}\cite{CCD}\cite{CCD2}\cite{CWZ}\cite{SD} etc..
 \begin{figure}
\includegraphics[width=15.0true cm, height=9.00true cm]{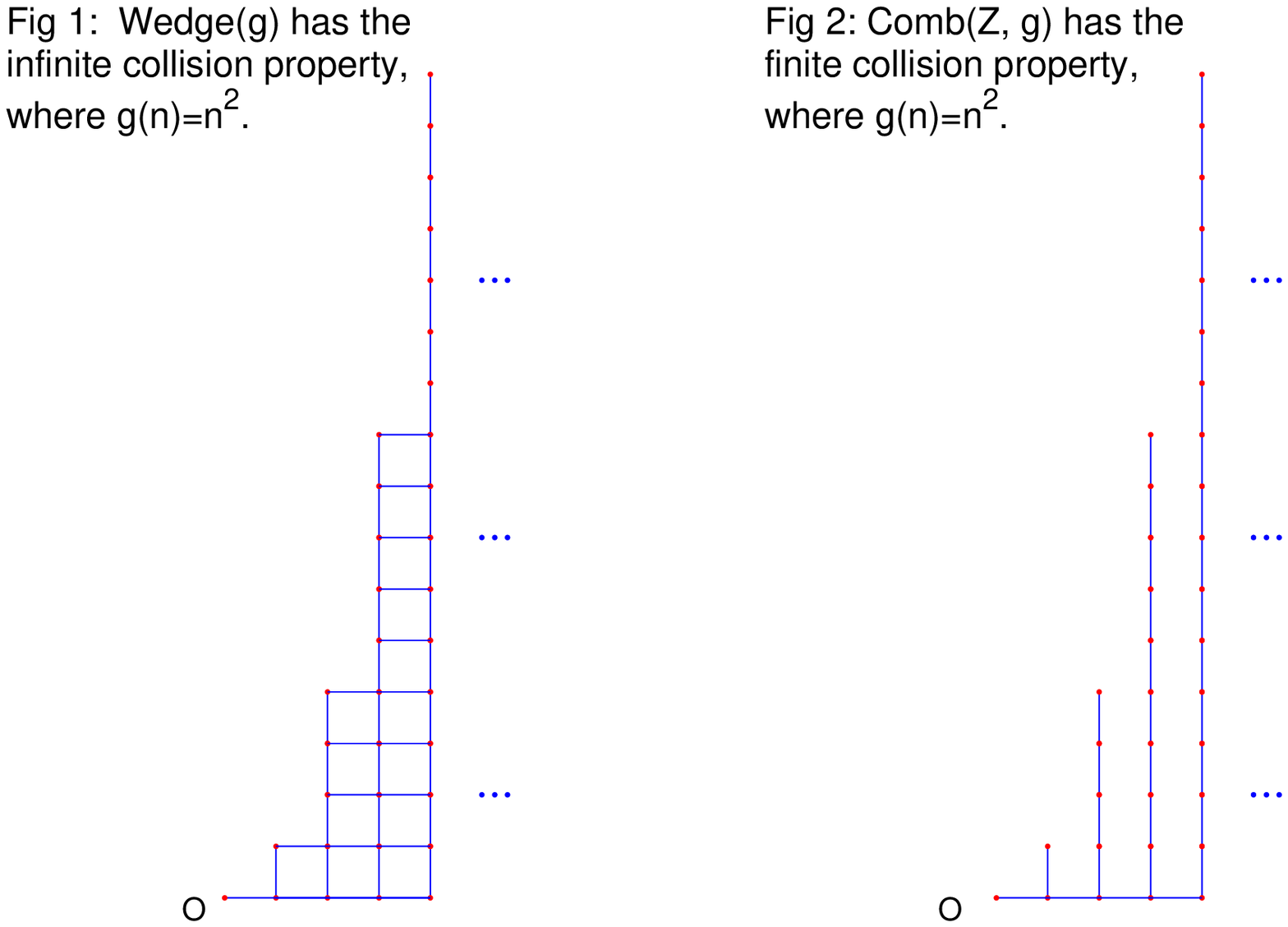}
\end{figure}
\begin{thm}\label{t:1.1}
Wedge$(f_1,\cdots,f_d)$ has the infinite collision property if and
only if Wedge$(f_1,\cdots,f_d)$ is recurrent.
\end{thm}
To understand the conditions better, it is worthwhile to compare a
wedge with a wedge comb. Wedge$(g)$ always has the infinite
collision property since any subgraph of $\ZZ^2$ is recurrent.
However, Comb($\ZZ, g$) may have the finite collision property
\cite{BP}\cite{KP}. Refer to Figure 1 and Figure 2.   It implies
that our theorem holds owing to the monotone property of the profile
$f_i(\cdot)$ of the wedge.


\section{A partition of vertex set $\VV$}
Obviously, the functions $h_1,\cdots, h_d$ defined in Section 1
satisfy that for each $1\le i\le d$ and each $n\ge 0$,
\begin{equation}\label{e:876}0\le h_i(n)\le f_i(n) {\rm~~~and~~~} h_i(n+1)-h_i(n)\in\{0
,1\}. \end{equation} We shall define a class of subsets
$\Delta_i(n)$ and $\partial_n$ through  these functions. We shall
show later that $\{\partial_n: n\ge 0\}$ is a partition of
$\VV$. 
 For each $1\le i\le
d+1$, let
$$
\Delta_i(0)=\{(0,\cdots, 0)\}\in \ZZ^{d+1}.
$$
 Fix $n\ge 1$, let
$$
\Delta_{d+1}(n)=\{(x_1,\cdots, x_d, n)\in \ZZ^{d+1}: 0\le x_i\le
h_i(n), 1\le i\le d\}.
$$
Then $\Delta_{d+1}(n)$ is a subset of $\VV$.  Fix  $1\le i\le d$. If
$h_i(n)=h_i(n-1)+1$ then let
$$
\Delta_i(n)=\{(x_1,\cdots, x_d,x_{d+1})\in\VV:  x_j\le h_j(n) {\rm
~for ~each~} 1\le j\le d,~ x_i=h_i(n),~ x_{d+1} \le n \}.
$$
Otherwise, if $h_i(n)=h_i(n-1)$ then let $\Delta_i(n)=\emptyset$.\\
For each $n\ge 0$ we set
$$
\partial_n=\bigcup_{i=1}^{d+1} \Delta_i(n).
$$
Finally, for each $x\in \RR^{d+1}$ and each $1\le i\le d+1$,  we
denote by $x_{i}$ the $i$-th coordinate of $x$. For each $x\in \VV$
and $1\le i\le d$, we set
$$
p_i(x)=\min\{m: h_i(m)\ge x_i\}.
$$
By (\ref{e:876})
$$
h_i(p_i(x))=x_i.
$$
For each $x\in \VV$, set $$u(x)=\max\{x_{d+1}, p_1(x),\cdots,
p_d(x)\}.$$ Then we have the following lemma.
\begin{lemma} For each pair of $m\ge 0$ and $x\in \VV$,  vertex $x\in \partial_m$
if and only if $u(x)=m$.
\end{lemma}
\pf   Fix $x=(x_1,\cdots, x_d,n)\in\VV$. For conciseness, we write
$p_i$ instead of $p_i(x)$. First   we shall prove the statement that
if $u(x)=m$ then $x\in
\partial_m$. Set $$S=\{ i: 1\le i\le d, x_i> h_i(n)\}.$$ We consider two
cases $S=\emptyset$ and $S\not=\emptyset$.

Case I:  $S=\emptyset$.  Then for each $1\le i\le d$,
$$
x_i\le h_i(n).
$$
As a result,
$$x\in \Delta_{d+1}(n)\subset \partial_n.$$
Since $ h_i(p_i)=x_i$,
$$
h_i(p_i)\le h_i(n).
$$
By the definition of $p_i(\cdot)$,
$$
p_i\le n.
$$
 Therefore, $u(x)=n$ as claimed above.

Case II: $S\not= \emptyset$. Fix $j\in S$ which satisfies that for
all $l\in S$,
\begin{equation}\label{e:2035} p_l\le p_j.
\end{equation}
We shall show that $u(x)=p_j$ and $x\in
\partial_{p_j}$.
Since $j\in S$,
$$
h_j(p_j)=x_j>h_j(n).
$$
It implies that
\begin{equation}\label{e:205}
n< p_j.
\end{equation}
 Furthermore, for each $l\in \{1,\cdots, d\}\setminus S$
\begin{equation}\label{e:207}
h_l(p_l)=x_l\le h_l(n)\le h_l(p_j).
\end{equation}
As a result of that
\begin{equation}\label{e:2075}
p_l\le p_j.
\end{equation}
Owing to (\ref{e:2035}), (\ref{e:205}) and (\ref{e:2075}),
$$
u(x)=p_j.
$$
On the other hand,  by the definition of $p_j(\cdot)$ there has

$~~~~~~~~~~~~~~~~~~~~~~~~~$either $~~~p_j=0~~~~$ or $~~~~h_j(p_j-1)< h_j(p_j)$.\\
However, there always have
\begin{equation}\label{e:206}\Delta_j(p_j)=\{(y_1,\cdots, y_d,y_{d+1})\in\VV: y_l\le h_l(p_j)
{\rm ~for ~each~} 1\le l\le d,~y_j=h_j(p_j), ~ y_{d+1} \le p_j \}.
\end{equation}
By (\ref{e:2035}), for each $l\in S$
\begin{equation}\label{e:204}
x_l= h_l(p_l)\le h_l(p_j).
\end{equation}
 By (\ref{e:207}), (\ref{e:206}) and (\ref{e:204}), we have
that
$$
x\in \Delta_j(p_j)\subset \partial_{p_j}.
$$
Such we have proved the first statement for both cases. \\

 Next  we shall show that $\partial_0,\partial_1,\cdots$ are disjoined. Fix $n>m\ge 0$.
 Since that for any $x\in
\Delta_{d+1}(n)$ and any $y\in \partial_m$, $$ x_{d+1}=n>m\ge
y_{d+1}.
$$ So,
\begin{equation}\label{e:200} \partial_m\cap \Delta_{d+1}(n)=\emptyset.\end{equation} Fix $1\le i\le d$ and $ 1\le j\le d$. We will show that $
\Delta_{i}(m)\cap \Delta_{j}(n)= \emptyset$. Otherwise, suppose  $
\Delta_{i}(m)\cap \Delta_{j}(n)\not= \emptyset$. Then
$$
\Delta_i(m)=\{(x_1,\cdots, x_d,x_{d+1})\in\VV:  x_l\le h_l(m) {\rm
~for ~each~} 1\le l\le d,~ x_i=h_i(m),~ x_{d+1} \le m \},
$$
$$
\Delta_j(n)=\{(x_1,\cdots, x_d,x_{d+1})\in\VV:  x_l\le h_l(n) {\rm
~for ~each~} 1\le l\le d,~ x_j=h_j(n),~ x_{d+1} \le n \}.
$$
And then
$$
h_j(n)=h_j(n-1)+1.
$$
 Furthermore, since $
\Delta_{i}(m)\cap \Delta_{j}(n)\not= \emptyset$ there exists $z\in
\Delta_{i}(m)\cap \Delta_{j}(n)$. Then
$$
z_j=h_j(n) {\rm~and~} z_l\le \min\{h_l(m), h_l(n)\} {\rm~ for ~each~
}1\le l\le d.
$$
Hence,
\begin{equation}\label{e:512}
h_j(n)\le h_j(m).
\end{equation}
On the other hand,  since $h_j(\cdot)$ is an increasing function and
$n>m$,
$$
h_j(n-1)\ge h_j(m).
$$
It deduces that $$ h_j(n)=h_j(n-1)+1\ge h_j(m)+1>h_j(m).
$$ This contradict  (\ref{e:512}). Therefore,
\begin{equation}\label{e:202}\Delta_{i}(m)\cap
\Delta_{j}(n)=\emptyset.\end{equation}
 Similarly, we can
prove that
\begin{equation}\label{e:201}
\Delta_i(n)\cap \Delta_{d+1}(m)= \emptyset.
\end{equation}
 Taking (\ref{e:200}), (\ref{e:202}) and
(\ref{e:201}) together, we get that $\partial_n$ and $\partial_m$
are disjoined. We have finished the proof of the lemma.\qed
\\

The next lemma shows that the neighbor of $\partial_{n}$ are
$\partial_{n-1}$ and $\partial_n$  for each $n\ge 1$. It implies
that $\partial_n$ is a cutset of the graph Wedge($f_1,\cdots, f_d$).
We write
$$e_i=(0,\cdots,0,1,0,\cdots,0)$$ for the $i$-th unit vector of
$\RR^{d+1}.$

\begin{lemma}\label{l:532} Let $x\in \VV$ and $1\le i\le d+1$. If $x+e_i\in \VV$
then
$$
u(x+e_i)-u(x)=0  {\rm~or~}1.
$$
\end{lemma}
\pf  Fix  $x\in \VV$.  Obviously for each $1\le i\le d+1$ and $1\le
l\le d+1$ with $i\not=l$, if $x+e_l\in\VV$ then
$$p_i(x+e_{l})=p_i(x).$$
First we consider the easy case $i=d+1$. Obviously,
$x+e_{d+1}\in\VV$. Hence
\begin{align*}
u(x+e_{d+1})-u(x)=&\max\{x_{d+1}+1, p_1(x+e_{d+1}),\cdots,
p_d(x+e_{d+1})\}-\max\{x_{d+1}, p_1(x),\cdots, p_d(x)\}\\
 =& \max\{x_{d+1}+1, p_1(x),\cdots,
p_d(x)\}-\max\{x_{d+1}, p_1(x),\cdots, p_d(x)\}\\ =& ~0 {\rm~or ~}1.
\end{align*}
Next we consider the case $1\le i\le d$. Fix  $x\in \VV$ and
$x+e_i\in \VV$. \\
If $f_i(p_i(x)+1)\ge x_i+1$, then
$$
f_i(p_i(x)+1)\ge  x_i+1=h_i(p_i(x))+1.
$$
Hence $$h_i(p_i(x)+1)=h_i(p_i(x))+1=x_i+1.$$ Such
$$
p_i(x+e_i)=p_i(x)+1.
$$
Similarly we have
\begin{align*}
&u(x+e_{i})-u(x)\\
=&\max\{x_{d+1}, p_1(x+e_{i}),\cdots,
p_d(x+e_{i})\}-\max\{x_{d+1}, p_1(x),\cdots, p_d(x)\}\\
 =& \max\{x_{d+1}, p_1(x),\cdots,p_{i-1}(x), p_i(x)+1, p_{i+1}(x),
 \cdots,
p_d(x)\}-\max\{x_{d+1}, p_1(x),\cdots, p_d(x)\}\\ =& ~0 {\rm~or ~}1.
\end{align*}
Otherwise,    $f_i(p_i(x)+1)< x_i+1$. Let
$$\eta_i=\min\{m: f_i(m)\ge x_i+1\}.   $$ Then
$$\eta_i>
p_i(x)+1.$$ Furthermore,
$$
h_i(\eta_i-1)\ge h_i(p_i(x))=x_i.
$$
On the other hand
$$
h_i(\eta_i-1)\le f_i(\eta_i-1)<x_i+1.
$$
Since $h_i(\cdot)$ is integer valued,
$$
h_i(\eta_i-1) = x_i.
$$
As a result,
$$
f_i(\eta_i)\ge x_i+1= h_i(\eta_i-1)+1.
$$
Hence
$$
h_i(\eta_i)=h_i(\eta_i-1)+1=x_i+1.
$$
Therefore,
\begin{equation}\label{e:441}
p_i(x+e_i)\le \eta_i.
\end{equation}
Since $x+e_i\in\VV$,
$$
f_i(x_{d+1})\ge x_i+1.
$$
and then $$ \eta_i\le x_{d+1}.
$$
By (\ref{e:441}), $$ p_i(x+e_i)\le x_{d+1}.
$$
So that,
\begin{align*}
&u(x+e_{i})-u(x)\\
=&\max\{x_{d+1}, p_1(x+e_{i}),\cdots,
p_d(x+e_{i})\}-\max\{x_{d+1}, p_1(x),\cdots, p_d(x)\}\\
 \le & \max\{x_{d+1}, p_1(x),\cdots,p_{i-1}(x), x_{d+1},  p_{i+1}(x),
 \cdots,
p_d(x)\}-\max\{x_{d+1}, p_1(x),\cdots, p_d(x)\}\le 0.
\end{align*}
By the increasing property of $u(\cdot)$, we get that
$$u(x+e_i)-u(x)=0. $$
 \qed
\\

At the end of this section, we shall estimate the cardinality of
$\partial_n$.
\begin{lemma}\label{l:2.222} For each $n\ge 0$,  $$\prod_{i=1}^d (h_i(n)+1)\le |\partial_n|\le
(d+1) \prod_{i=1}^d (h_i(n)+1) .$$
\end{lemma}
\pf For each $n\ge 0$
$$
|\partial_n|\ge |\Delta_{d+1}(n)|=\prod_{i=1}^d (h_i(n)+1),
$$
since $ \Delta_{d+1}(n) \subseteq \partial_n$.

Fix $n\ge 1$ and  $1\le i\le d$. Without making confusion, we set
$$p_i=p_i(n)=\min\{m: h_i(m)=n\}.$$
Then
\begin{align*}\Delta_i(p_i)=&\{(x_1,\cdots, x_d,x_{d+1})\in\VV: x_l\le h_l(p_i)
{\rm ~for ~each~} 1\le l\le d,~ x_i=n,~ x_{d+1} \le p_i \}.
\end{align*}
As we have known that if $x\in \VV$ with $x_i=n$ then $
f_i(x_{d+1})\ge n.$
 Let $$k=\min\{u\in \ZZ^+: f_i(u)\ge
n\}.$$ Then
\begin{align*}\Delta_i(p_i)
= &\{(x_1,\cdots, x_d,x_{d+1})\in\VV: 0\le x_l\le h_l(p_i) {\rm ~for
~each~} 1\le l\le d,~ x_i=n,~ k\le x_{d+1} \le p_i\}\\
\subseteq & \{(x_1,\cdots, x_d,x_{d+1})\in\ZZ^{d+1}: 0\le x_l\le
h_l(p_i) {\rm ~for ~each~} 1\le l\le d,~ x_i=n,~ k\le x_{d+1} \le
p_i\}.
\end{align*}
Therefore,
$$
|\Delta_i(p_i)|\le\frac{p_i-k+1}{h_i(p_i)+1} \prod_{l=1}^{d}
(h_l(p_l)+1).
$$
If $k\le \eta<p_i$, then $$h_i(\eta)+1\le h_i(p_i-1)+1=h_i(p_i)=n\le
f_i(k)\le f_i(\eta).$$ And then
$$
h_i(\eta)=h_i(\eta-1)+1.
$$
Therefore,
$$
h_i(p_i)-h_i(k)=p_i-k.
$$
Such  $$|\Delta_i(p_i)|\le \frac{h_i(p_i)-h_i(k)+1}{h_i(p_i)+1}
\prod_{l=1}^{d} (h_l(p_i)+1)\le \prod_{l=1}^{d} (h_l(p_i)+1).$$ So
that for any $m\ge 0$ , if $m\in\{p_i(n):n\ge 1\}$, then
\begin{equation}\label{e:3214112}
|\Delta_i(m)|\le \prod_{l=1}^{d} (h_l(m)+1).
\end{equation}

Obviously, (\ref{e:3214112}) is true for  $m=0$ since
$\Delta_i(0)=\{(0,\cdots,0)\}$.
 Notice that $p_i(0)=0$  and the fact that
 if $m\in \ZZ\backslash \{p_i(n):n\ge 0\}$ then
$ \Delta_i(m)=\emptyset.$ Therefore,  (\ref{e:3214112}) are true for
all $m\ge 0$. Finally, for any $m\ge 0$
$$
|\partial_m|\le \sum_{i=1}^{d+1}
|\Delta_i(m)|\le\sum_{i=1}^{d+1}\prod_{l=1}^{d} (h_l(m)+1)\le
(d+1)\prod_{i=1}^{d} (h_i(m)+1).
$$
 We have completed the proof
of the lemma.\qed

\section{Proof of Theorem \ref{t:1.2}}
We shall use the notation of electric network. Every edge of
Wedge$(f_1,\cdots,f_d)$ is assigned a unit conductance. So that, we
get an electric network.
 For sets $A,B\subset \VV$ with
$A\cap B=\emptyset$, denote by $\mathcal{R}(A\leftrightarrow B)$ the
effective resistance between $A$ and $B$ in the electric network.
For simplicity, we label $O$ as the origin of $\ZZ^{d+1}$ and  set
$$
\VV_r=\bigcup_{n=0}^r \partial_r
$$
for each $r\ge 1$. Then
 we have the following lemma.
\begin{lemma}\label{l:1.1} For each $r\ge 1$
$$
\mathcal{R}(O\leftrightarrow \partial_r)\ge
\frac{1}{2(d+1)^2}\sum_{n=0}^{r-1}\prod_{i=1}^d
\frac{1}{h_i(n)+1}.$$
\end{lemma}
\pf Notice that $\partial_0=\{O\}$.  By Lemma \ref{l:532}, for each
$n\ge 1$ the neighbor of $\partial_n$ are $\partial_{n-1}$ and
$\partial_{n+1}$
 in Wedge($f_1,\cdots, f_d$). So that $\partial_n$
is a cutset which separates $O$ from $\partial_{n+s}$. The rest
proof is easy and one can refer to \cite{Peres}. Fix $r$. The
effective resistance from $O$ to $\partial_r$ in $(\VV,\EE)$ is
equal to that in  its subgraph with vertex set $\VV_r$. We short
together all the vertices in $\partial_n$ for each $0\le n\le r$.
And replace the edges between $\partial_{n}$ and $\partial_{n+1}$ by
a single edge of resistance $\frac{1}{b_n}$, where $b_n$ is the
number of edges connect $\partial_{n}$ with $\partial_{n+1}$. This
new network is a series network with the same effective resistance
from $O$ to $\partial_r$. Thus, Rayleigh's monotonicity law shows
that the effective resistance from $O$ to $\partial_r$ in $\VV_r$ is
at least $\sum_{n=0}^{r-1} \frac{1}{b_n }$. By  Lemma \ref{l:2.222}
and the fact that every vertex of Wedge($f_1,\cdots, f_d$) has at
most $2(d+1)$ neighbor,
 $$
\mathcal{R}(O\leftrightarrow \partial_r)\ge \sum_{n=0}^{r-1}
\frac{1}{b_n }\ge \frac{1}{2(d+1)}\sum_{n=0}^{r-1}
\frac{1}{|\partial_n| }\ge \frac{1}{2(d+1)^2}\sum_{n=0}^{r-1}
\prod_{i=1}^d\frac{1}{h_i(n)+1}.
$$  \qed
\\

On the other hand we can estimate the upper bound of
$\mathcal{R}(x\leftrightarrow \partial_r)$.

\begin{lemma}\label{l:5.12} There exists $C_d>0$ which depends only on $d$ such that  for any $r\ge
1$ and any  $x\in \VV_{r-1}$,
$$
\mathcal{R}(x\leftrightarrow
\partial_r)\le  C_d\sum_{n=0}^{r-1}\prod_{i=1}^d\frac{1}{h_i(n)+1}.$$
\end{lemma}
\pf Outline of the proof. We shall  construct $2d$ functions $g_{\pm
i}(\cdot)$ first. These functions will help us  to find a subset
$\VV_x$ which satisfies that $x\in\VV_x\subseteq \VV_r$. Such
$\mathcal{R}_{\VV_x}(x\leftrightarrow \Delta_{d+1}(r)\cap\VV_x)$,
the resistance between $x$ and $\Delta_{d+1}(r)\cap\VV_x$ in the
subgraph with vertex set $\VV_x$, is greater than
$\mathcal{R}(x\leftrightarrow
\partial_r)$.
 Furthermore, we show the relation between  $\VV_x$ and
Wedge($h_1,\cdots, h_d$). As  known from Lyons\cite{TL}, the related
resistance   in Wedge($h_1,\cdots, h_d$) can be gotten. So do
$\mathcal{R}_{\VV_x}(x\leftrightarrow
\Delta_{d+1}(r)\cap\VV_x)$.\\

Fix $x=(x_1,\cdots,x_{d}, s)\in \VV_{r-1}$. We shall construct $2d$
nonnegative integer valued functions on $\ZZ^+$.
 Fix $1\le i\le
d$.  First set $$g_{\pm i}(0)=x_i.$$ Suppose  that  the definition
of $g_{\pm i}(n)$ is known, we define $g_{\pm
i}(n+1)$ in three cases.\\
 (1) If
$h_i(n+1)=h_i(n)$, then we set $g_{\pm i}(n)=g_{\pm i}(n+1)$.\\
(2)
 If
$h_i(n+1)=h_i(n)+1$ and if $g_{-i}(n)=0$, then we set
$g_{-i}(n+1)=0$ and $g_{i}(n+1)=g_i(n)+1$.\\
(3)
 Otherwise, if $h_i(n+1)=h_i(n)+1$ and if $g_{-i}(n)>0$, then we set
$g_{-i}(n+1)=g_{-i}(n)-1$ and $g_{i}(n+1)=g_i(n)$.\\
We  say that these functions $g_{\pm i}(n)$ has the properties
(a),(b) and (c). Where
\begin{align*}
(a)&: g_{i}(n+1)-g_{i}(n)\in \{0,1\} {\rm ~and~
}g_{-i}(n+1)-g_{-i}(n)\in \{0,-1\} {\rm~for~each~} n\ge 0;\\
 (b)&: g_{i}(n)-g_{-i}(n)=h_{i}(n) {\rm~for~each~} n\ge 0;\\
(c)&: 0\le g_{-i}(n)\le g_{i}(n)\le \min\{ f_{i}(n+s),
h_i(r)\}{\rm~for~each~} 0\le n\le r-s.
\end{align*}
 Obviously, (a) are
true for all $n\ge 0$. Next we shall prove (b) by induction to $n$.
It is true for  $n=0$ since $h_i(0)=0$. Suppose (b) is true for
$n=m$ and we shall check $n=m+1$. In any case of (1),(2) and (3),
there has
$$
h_{i}(m+1)-h_i(m)=[g_i(m+1)-g_i(m)]-[g_{-i}(m+1)-g_{-i}(m)].
$$
By the assumption that (b) is true for $n=m$, we can get that (b) is
still true for $n=m+1$. Such (b) is true for any $n\ge 0$. Again we
prove (c) by induction. Owing to $x\in \VV_{r-1}$ and $x_{d+1}=s$,
$$
0\le x_i\le h_i(x_{d+1})=h_i(s)\le \min\{h_i(r), f_i(s)\}.
$$
So  (c) is true  for $n=0$.  Suppose (c) is true for $n=m<r-s$ and
we
shall check $n=m+1$. \\
If (1) is true for $n=m+1$, then by the assumption that (c) is true
for $n=m$ and the  monotone property of  $f_i(\cdot)$, we have
(c)  for $n=m+1$. \\
If (2) is true for $n=m+1$, then what we need to care is only
$g_i(n+1)$. However, by the result (b) we have proved
$$
g_i(n+1)=h_i(n+1)+g_{-i}(n+1)=h_i(n+1)\le f_i(n+1)\le f_i(s+n+1).
$$
Furthermore, since $n<r-s$, $$h_i(n+1)\le  h_i(r).
$$
Therefore  (c) is true for $n=m+1$. \\
 If (3) is true for $n=m+1$,
then what we need to care is only $g_{-i}(n+1)$. But by the
condition that $g_i(n)>0$, we have
$$
g_{-i}(n+1)=g_{-i}(n)-1\ge 0.
$$
Hence (c) is true, too. Therefore, in any case (c) is true for
$n=m+1$ with $n<r-s$.\\

As a result, we can define vertex set $\VV_x$ and edge set $\EE_x$.
Let
$$
{\VV}_x=\{(u_1,\cdots,u_d, n+s)\in \ZZ^{d+1}: 0\le n\le r-s,~
g_{-i}(n)\le u_i\le g_{i}(n) {\rm~for~each~} 1\le i\le d
 \}.
$$
Let $${\EE}_x=\{[u,v]\in\EE: u,v\in {\VV}_x\}.$$ The definition does
not make confusion of $\VV_x$ and $\VV_n$ since $x$ is a vector.
 By (c), $$x\in{\VV}_x\subseteq\VV_r.$$
Hence graph $({\VV}_x, {\EE}_x)$ is a subgraph of
Wedge$(f_{1},\cdots, f_{d})$. Notice that $$
 \partial_r\cap
{\VV}_x\supseteq\Delta_{d+1}(r)\cap\VV_x.
$$
(Actually   $\partial_r\cap {\VV}_x=\Delta_{d+1}(r)\cap\VV_x$, but
we omit the proof here since it is irrelevant to our main result.)
By the Rayleigh's monotonicity law, the effective resistance between
$x$ and $\Delta_{d+1}(r)\cap\VV_x$ in the subgraph
 is greater than that in the old graph. That is,
\begin{equation}\label{e:3.22}
\mathcal{R}(x\leftrightarrow
\partial_r)\le\mathcal{R}_{{\VV}_x}(x\leftrightarrow
\Delta_{d+1}(r)\cap {\VV}_x).
\end{equation}
So that we need only to estimate the upper bound of
$\mathcal{R}_{{\VV}_x}(x\leftrightarrow \Delta_{d+1}(r)\cap
{\VV}_x)$.\\

We shall show the relation between  $(\VV_x,\EE_x)$ and
Wedge($h_1,\cdots,h_d$). Let
$$
\mathbb{H}=\{(x_1,\cdots, x_{d}, n)\in\ZZ^{d+1}: 0\le x_i\le h_i(n)
{\rm ~for ~each~ }1\le i\le d, 0\le n\le r-s\}.
$$
Obviously, $ \mathbb{H}$ is a subset of  vertices of
Wedge($h_1,\cdots,h_d$).  By the construction of $g_{\pm i}(\cdot)$,
one can easily check that there has for
 each $n\ge 1$\\
 $~~~~~~~~~~~~~~~~~~~~~~~$either
$~~~g_{-i}(n)=g_{-i}(n-1)~~~~$ or $~~~~g_{i}(n)=g_{i}(n-1)$.
\\
 So we can define
$$
L_i(n)=\min\{g_{si}(n): g_{si}(n)=g_{si}(n-1), s\in\{-1,1\}\}.
$$
Let $\Gamma(x)=O$.  For each $(u_1,\cdots, u_d, n+s)\in {\VV}_x$
with $n\ge 1$, let
$$
\Gamma(u_1,\cdots, u_d, n+s)=(|u_1-L_1(n)|,\cdots, |u_d-L_d(n)|, n).
$$
By (b), $\Gamma$ is  a bijection function from ${\VV}_x$
 to $\mathbb{H}$. Obviously, $[u,v]\in\EE_x$ if and only
if $[\Gamma(u), \Gamma(v)]$ is an edge of  Wedge($h_1,\cdots,h_d$)
 for each pair of
  $u$ and $v$ with $u_{d+1}=v_{d+1}$. Moreover, for any $u\in\VV_x$ we have that  $u-e_{d+1}\in \VV_x$ if and only if $\Gamma(u)-e_{d+1}\in\HH$. \\

Since $h_i(\cdot)$ increases at most one at each step, we can use
the result of Lyons\cite{TL}. That is,  there exists a unit flow
$\mathbf{w}$ from $O$ to $\Delta_{d+1}(r-s)$ in the  subgraph of
Wedge($h_1,\cdots,h_d$) with vertex set $\HH$, such that for each
$u\in \HH$ with $u_{d+1}=n<r-s$,
\begin{equation}\label{e:01022}
\mathbf{w}(u,u+e_{d+1})=\prod_{i=1}^d\frac{1}{h_i(n)+1},
\end{equation}
and the energy of $\mathbf{w}$ has upper bound
\begin{equation}\label{e:21201}
\mathcal{E}(\mathbf{w})\le C_d
\sum_{n=0}^{r-s-1}\prod_{i=1}^d\frac{1}{h_i(n)+1},
\end{equation}
where $C_d<\infty$ and depends only on $d$. Let $\mathbf{{w}}_x$ be
a function on ${\EE}_x$ and satisfies that for each $[u,v]\in
{\EE}_x$ with $u_{d+1}=v_{d+1}$,
\begin{align*}
\mathbf{{w}}_x(u,v)=\mathbf{{w}}(\Gamma(u),\Gamma(v)).
\end{align*}
and for each $u\in \VV_{r-1}$ with $u_{d+1}=n$, let $$
\mathbf{{w}}_x(u,u+e_{d+1})=\prod_{i=1}^d\frac{1}{h_i(n)+1}.
$$
Directly calculate
\begin{align*}
&~\sum_{v:[u,v]\in\EE_x}\mathbf{{w}}_x(u,v)\\
&=\mathbf{{w}}_x(u,u+e_{d+1})+\mathbf{{w}}_x(u,u-e_{d+1})1_{\{u-e_{d+1}\in\VV_x\}}+\sum_{v:[u,v]\in\EE_x,
u_{d+1}=v_{d+1}}\mathbf{{w}}_x(u,v)\\
&=\prod_{i=1}^d\frac{1}{h_i(n)+1}-\prod_{i=1}^d\frac{1}{h_i(n-1)+1}1_{\{u-e_{d+1}\in\VV_x\}}+\sum_{v:[u,v]\in\EE_x,
u_{d+1}=v_{d+1}}\mathbf{{w}}(\Gamma(u),\Gamma(v))\\
&=\mathbf{{w}}(\Gamma(u),\Gamma(u)+e_{d+1})+\mathbf{{w}}(\Gamma(u),\Gamma(u)-e_{d+1})1_{\{\Gamma(u)-e_{d+1}\in\HH\}}+\sum_{
z\in \HH:\|u-z\|_1=1,
u_{d+1}=z_{d+1}}\mathbf{{w}}(\Gamma(u),z)\\
&=\sum_{ z\in \HH:\|u-z\|_1=1}\mathbf{{w}}(\Gamma(u),z).
\end{align*}
Together with the fact that  $\mathbf{{w}}$ is a unit flow, we get
that $\mathbf{{w}}_x$ is a unit flow from $x$ to
$\Delta_{d+1}(r)\cap \VV_x$ in graph $({\VV}_x, {\EE}_x)$. Obviously
\begin{equation}\label{e:01012}
\mathcal{E}(\mathbf{{w}}_x)=\mathcal{E}(\mathbf{{w}}).
\end{equation}
Together   (\ref{e:3.22}), (\ref{e:21201}) and (\ref{e:01012}), we
have
$$
\mathcal{R}(x\leftrightarrow\partial_r)\le
\mathcal{R}_{{\VV}_x}(x\leftrightarrow \Delta_{d+1}(r)\cap
\VV_x)\le\mathcal{E}(\mathbf{{w}}_x)=\mathcal{E}(\mathbf{{w}})\le
C_d\sum_{n=0}^{r-1}\prod_{i=1}^d\frac{1}{h_i(n)+1}.
$$
\qed
\\

{\it Proof of Theorem \ref{t:1.2}.} As it is well known, a connect
graph with local finite degree is recurrent if and only if the
resistance from  any one vertex  to the infinity in the graph is
infinite (Refer to \cite{Peres},  Proposition 9.1). Together with
Lemmas \ref{l:1.1} and \ref{l:5.12}, we have the desired result.\qed

\section{Proof of Theorem \ref{t:1.1}}

\begin{lemma}\label{l:1.3} Let $G$ be a graph of bounded degrees with a
distinguished vertex $o$ and suppose that there exists a sequence of
sets $(B_r)_r$ growing with $r$ and satisfying
$$
g_{B_r}(o,o)\rightarrow\infty {\rm~as~}r\rightarrow\infty~{\rm
and~}g_{B_r}(x,x)\le Cg_{B_r}(o,o),~~\forall x\in G,
$$
for a uniform constant $C>0$. Here, $g_{B}(\cdot,\cdot)$ is the
green function of the simple random walk on $G$ killed when it exits
$B$. Then the graph $G$ has the infinite collision property.
\end{lemma}
\pf Refer to \cite{BP}. \qed\\

{\it Proof of Theorem \ref{t:1.1}.} First suppose
Wedge($f_1,\cdots,f_d)$ is not a recurrent graph. Then
Wedge($f_1,\cdots,f_d)$  is a transient graph. It implies that
$g_{\VV}(O,O)$, the expected number of returning to $O$, is finite.
One can easily get that the expected number of collisions between
two independent simple random walks starting from $O$ is less than
$2(d+1)g_{\VV}(O,O)$. So that, almost surely the number of
collisions  is finite. Hence, Wedge($f_1,\cdots,f_d)$ has the finite
collision property.

On the other hand, suppose Wedge($f_1,\cdots,f_d)$ is recurrent. By
Theorem \ref{t:1.2} we have (\ref{e:1.1}). Furthermore,  by Lemma
\ref{l:1.1}
$$
\lim_{r\rightarrow\infty}\mathcal{R}(O\leftrightarrow
\partial_{r})\ge \lim_{r\rightarrow\infty}\frac{1}{2(d+1)^2}\sum_{n=0}^{r-1}\prod_{i=1}^d\frac{1}{f_i(n)+1}=\infty.
$$
 As it is known to
all (refer to \cite{BP}) that for each $r\ge 1$
$$
\mathcal{R}(O\leftrightarrow \partial_{r+1})=g_{\VV_{r}}(O,O).
$$
So $\lim_{r\rightarrow\infty}g_{\VV_r}(O,O)=\infty$. By Lemmas
\ref{l:1.1} and \ref{l:5.12}, for all $r\ge1$ and $x\in {\rm
Wedge}(f_1,\cdots,f_d)$
$$
g_{\VV_r}(x,x)\le 2(d+1)^2C_d~g_{\VV_r}(O,O).
$$
By Lemma \ref{l:1.3}, Wedge($f_1,\cdots,f_d)$  has the infinite
collision property.\qed

\end{document}